\documentclass{book}
\usepackage[contrib,lang=british %french
]{ems-book} %% change to `american' if you use American English

\usepackage{tikz-cd}% Diagram

\theoremstyle{plain}
\newtheorem{theorem}{Theorem}[chapter]

\theoremstyle{definition}

\newtheorem{remark}[theorem]{Remark}
\newtheorem{definition}[theorem]{Definition}

\newcommand{\R}{\mathbb{R}}
\newcommand{\N}{\mathbb{N}}

\newcommand{\CC}{\mathbb{C}}

\begin{document}
\mainmatter

%------
% Insert the title of your paper and (if necessary)
% a short title for the running head.
%------
\title{An equisingular heritage of Bernard Teissier }
\titlemark{An equisingular heritage of Bernard Teissier}

%------
% Insert full names of the authors.
% Add further authors as follows:
%  \emsauthor{2}{}{}
%  \emsauthor{3}{}{}
% etc.
% Abbreviate first names for the running head.
%------
\emsauthor{1}{Georges Comte}{G.~Comte}

%------
% Use \authormark if the list of authors is too
% long for the running head: \authormark{A.~Doe et al.}
%------

%------
% Add one \emsaffil and one \email for each author.
% NOTE: The address does NOT appear in the paper.
% It will probably be printed in an appendix.
%------
\emsaffil{1}{
Laboratoire de Math\'ematiques de l'Universit\'e de Savoie Mont~Blanc, UMR CNRS 5127,
B\^atiment Chablais, Campus scientifique,
73376 Le Bourget-du-Lac cedex, France
\email{georges.comte@univ-smb.fr}
\url{https://georgescomte.perso.math.cnrs.fr/}
}

%------
% Add MSC 2020 codes according to www.ams.org/msc/msc2020.html.
% Secondary codes (in square brackets) are optional.
%------
\classification%[XXxXX]
{14B05,14B07, 32S60}

%------
% Add a list of keywords.
%------

%------
% Optional: dedication
%------
\chapterdedication{In honor of Bernard Teissier's 80 birthday}

%------
% Insert your abstract.
%------
\begin{abstract}

We give a brief and partial overview of Bernard Teissier's work in complex equisingularity theory, and a perspective on its legacy; in particular, we focus on the development of the theory in the real and the non-Archimedean contexts. Our aim is not to go into technical details, but, hopefully, rather to 
give a flavour of the forms taken by these developments, while providing enough definitions and references to give the reader access to old and new reference articles in the field.

\end{abstract}
%\tableofcontents
\makecontribtitle

\section{Introduction}
Exaggerating a bit, for many years, the easiest way to meet Bernard Teissier in the world was to climb the Ar\^etes de la Bruyère in the Cerces massif of the French Alps, since after the straight traverse of the ridge, rappelling
\footnote{de-escalation without a rope is possible from this point, but this detail is of no interest for the purposes we are pursuing here.} 
 from the last tower drop you off at the foot of the rustic shepherd's hut he lived in every summer with his wife.  Conversely, a direct approach to Bernard's hut from below, especially in foggy weather, was somewhat hazardous; all shepherd's huts are built to escape the intrusive sight of hikers. 
 
This is also how singularity theory works, a field of mathematics in which B. Teissier has played a major role since the 70s.  Indeed, understanding the geometry of a set at a singular point is nothing other than taking a roundabout route; the route of regularity.
The fact is that there is no primary, or intrinsic, mathematical definition of what a totally singular set should be, whereas there are several convincing definitions of what regular sets are.
%, and consequently convincing definitions of the breaks in regularity. 
%
The rules, and language games, from which mathematics proceeds, do not, by their very nature, allow the production of a discourse on what we call a singularity other than by defining various conditions of regularity, and then looking at the points where they fail.

%Therefore, in an apparently surprising way, presenting singularity theory a\-mounts first of all to reviewing relevant conditions of regularity, with a view then to violating them. 
The basic idea is the following: when a set $X$ no longer resembles, at a point $x\in X$, to what it resembles outside the point $x$ along a set $Y\subset X$ such that  $x\in \overline{Y}$, then the singular nature of this point emerges from the surrounding regularity. 

Equisingularity theory is the branch of mathematics where are investigated the notion of being alike for the germs $(X,y)_{y\in Y}$  along some set $Y$, and how to detect this propagation - and break - of similarity.
 Therefore Equisingularity theory is naturally based 
on the nested notions of stratifications and invariants. 
 On the one hand, stratification theory (see Section \ref{section. Equisingularity conditions}), provides the strata along which the set is alike, and on the other hand, invariants of singularities (see Section \ref{section. Numerical invariants}), detect, by their changes, jumps of strata.     
 Since there are as many notions of equisingularity as there are notions of  being alike for germs, that is as many notions of stratifications and singularity invariants, one of the main purposes of  Equisingularity theory is to compare all these notions of regularity by giving them a hierarchy.

There are excellent presentations of Equisingularity theory, from the beginning of the story, by B. Teissier himself (see \cite{Tei1975, Tei77}), to more recent and almost exhaustive ones (see among them \cite{Par2021, Bob2022}, where the first reference focuses on Zariski Equisingularity, and the second one more specifically on topological equisingularity). We focus here on some recent specific features of the development of the theory, such as the real and the non-Archimedean backgrounds. Those developments could be considered as complementary, or even sided developments of the theory, while we believe that they rather underline its flexibility and its fertility in, maybe, unexpected contexts.

%%%%%%%%%%%%%%%%%%%%%%%
\section{The logical-geometric structures}\label{subsection. The logicalgeometric structures}
%%%%%%%%%%%%%%%%%%%%%%%

As explained above, we aim at defining a singular point $x\in X$ by the degeneracy at $x$ of the regular geometry of $X$ along a regular set $Y\subset X$, such that  $x\in \overline{Y}$. In this goal the set $X$ has to belong to a class of sets for which having a regular geometry, whatever it means, is generic, or, in other words, a class of sets where the notion of dimension makes sense, and, moreover, where the codimension in $X$ of the set of non-regular points is positive. 
 In this direction, the more general framework is the framework of definable tame geometry, that is to say geometry where the tameness comes from the first order logical parsimony throughout which objects are defined. 

O. Zariski in \cite{Zar65a, Zar65b, Zar65c, Zar68, Zar71} initiated the theory in the context of algebraic or analytic varieties over algebraically closed ground fields of characteristic~$0$ (and sometimes  characteristic~$p\not=0$). In those classes, the algebraic, or the algebro-analytic, nature of objects provides the geometrical tameness, and gives in the same time a unified and powerful  language to translate this geometry. 

When one comes to subsets of $\R^n$, $n\in \N$, the largest convenient geometrical structures to consider are \emph{o-minimal structures}, encompassing the semialgebraic and the (globally) subanalytic structures, or for instance the structure $\R_{\mathrm{an,exp}}$, the expansion of the subanalytic structure by the exponential function. In such an o-minimal structure all the convenient geometrical operations and properties are at disposal: stability under finite product,  stability under finite union and finite intersection, projection, existence of dimension,  regular stratifications, smooth cell decomposition, finite number of connected components in family, and existence of an additive Euler-Poincaré characteristic etc. (see for instance \cite{vdDMil, vdD}).  

In Henselian valued fields, the non-archimedean counterpart of real o-minimal structures is given by \emph{Hensel minimality}, as defined in \cite{CluHalRid1, CluHalRid2}. The idea of Hensel minimality is to mimic the crucial finiteness axiom of o-minimal structures, that is, a definable  (in some given o-minimal structure on the real field $\R$) subset $X$ of $\R$ is a finite union of points and intervals, or put another way, there exist a finite set $C\subset \R$ and a map 
$\varepsilon:C\to \{-1,0,1\}$, depending on $X$, such that 
$$
x\in X \Longleftrightarrow \forall c\in C, \   \mathrm{sgn}(x-c)=\varepsilon(c).  
$$

Of course, in the non-archimedean context, a substitute of intervals, or of the $\mathrm{sgn}$ function has to be found. This is  done in \cite[Definition 1.2.3]{CluHalRid1} for the equi-characteristic $0$ case, and in \cite[Definition 2.2.1]{CluHalRid2} for the mixed characteristic case, in the following way (for sake of simplicity we only treat here the  equi-characteristic $0$ case).

Let $\mathcal{T}$ be a theory of non-trivially valued fields of equi-characteristic $0$ in a first order logical language $\mathrm{L}$ containing the language $\mathrm{L}_{\mathrm{val}} = \{ +, \cdot ,\mathcal{O} \}$ of valued fields (where the notation $\mathcal{O}$ is a notation for a predicate for the valuation ring). For $K$ a model of $\mathcal{T}$, we denote $\mathrm{val}_K$ its valuation, $\mathcal{O}_K$ its valuation ring and $\Gamma_K^\times$ its value group. Let $\lambda \le  1$ be an element of $\Gamma^\times_K$, and $I_\lambda = \{x \in  K;   \mathrm{val}_K(x) < \lambda \}$. Let us denote the quotient of multiplicative groups $K^\times / (1 + I_\lambda)$ by  $\mathrm{RV}^\times_\lambda$, and 
$$
\mathrm{rv}_\lambda : K \to  \mathrm{RV}_\lambda :=\mathrm{ RV}^\times_\lambda \cup  \{0\}
$$
the extension of the projection map $K^\times \to  \mathrm{RV}^\times_\lambda$, sending $0$ to $0$.
For $\ell\in \N\cup \{\omega\}$, the theory $\mathcal{T}$ is said
\emph{$\ell$-h-minimal} if for every model $K$ of $\mathcal{T}$, for every $\lambda \le 1$ in $\Gamma^\times_K$, for every set 
 $A \subset K$, for every set $A'\subset \mathrm{RV}_\lambda$ of cardinality $\le \ell$, and for every $(A\cup \mathrm{RV}_1 \cup A')$-definable set $X \subset K$, we have the existence of a $A$-definable set $C\subset K$ and a map $\varepsilon : C\to \mathrm{RV}_\lambda$ such that 
 $$
 x\in X \Longleftrightarrow \forall c\in C, \  \mathrm{rv}_\lambda(x-c)=\varepsilon(c). 
  $$
It is then shown in \cite{CluHalRid1, CluHalRid2} that sets definable in a structure of such an $\ell$-h-minimal theory (an Hensel minimal structure for short) have all the required tameness properties (that we don't give here).
%%%

Like the original approach of equisingularity by O. Zariski,  its immediate continuation by J. Briançon, J.-P. Henry, H. Hironaka, L\^e D\~ung Tr\'ang, I. Luengo, M. Merle, C. Sabbah, J.-P. Speder, B. Teissier and many others  (see \cite{BriSpe75, BriSpe76, BriGalGra81, HenMer81, HenMerSab84, Hir69, Lue87, Mer77, Spe75, Tei77a, LeTei81, Tei82, Var73}), is also mainly (complex) algebraic. This means that the definition of equisingular germs is given in purely algebraic terms, even if those authors have established deep connection between algebraic, differential, metric or topological equisingularity conditions. Of course, in the real and the non-Archimedean context, the algebraic approach has to be abandoned, and essentially remains the differential and metric approaches.    

In what follows we say that a set $X$ is a \emph{tame set} (or belonging to some \emph{tame geometric class}) if for some $n\in \N$, when $X\subset \CC^n$, $X$ is algebraic or analytic, or when $X\subset \R^n$, $X$ is definable in some o-minimal structure over the field $\R$, or when $X\subset K^n$, $X$ is definable in some Hensel minimal structure, with $K$ the valued field of this structure.  Notice that the identification of $\CC^n$ with $\R^{2n}$ turns any tame set of $\CC^n$ (at least locally in the complex analytic case) into a tame set of $\R^{2n}$.

%%%%%%%%%%%%%%%%%%%%%%%%%%%%%%%
\section{Equisingularity conditions }\label{section. Equisingularity conditions}
%%%%%%%%%%%%%%%%%%%%%%%%%%%%%%%

We now focus on various notions of being similar for two tame germs, or for a tame family of germs along a tame set.
\subsection{Topological equisingularity} 
The least demanding condition of similarity is the topological one; for $X,X'$ tame sets (complex or real, the non-archimedean case being here ruled out, since with no practical interest), we say that $X$ and $X'$ have the same (embedded) \emph{topological type at $x\in X$ and $x'\in X'$} if there exists a homeomorphism $(x,X,U) \to (x',X',U')$ (not necessarily belonging to the same tame class as $X$ and $X'$), where $U$ and $U'$ are open neighbourhoods of $x$ and $x'$ in the ambient space. 
And for $Y\subset X$ a smooth set (in the same tame class as $X$), we say that  $X$ is (locally) \emph{topologically equisingular along $Y$} when for any $y,y'\in Y$, $X$
has the same topological type at $y$ and $y'$. We know that there exists a finite decomposition of $X$ into smooth
\footnote{here smooth means $C^k$ for any chosen $k\in \N$, in a general o-minimal structure, whereas for the particular semi-algebraic or the subanalytic structures, one can take analytic smoothness.} 
pieces in the same class as $X$, along which a tame set is topologically equisingular.
This is for instance a direct consequence, in the o-minimal case, of the cell-decomposition theorem (see for instance \cite[Theorem 2.11, Chapter 3]{vdD}), and for the complex case, as we will see in subsection \ref{subsection. Whitney stratifications} below, we can refer more specifically to \cite[Chapiter VI-3]{Tei82}, where B. Teissier gave a construction of such a decomposition, based on numerical invariants. 
\subsection{Whitney stratifications}\label{subsection. Whitney stratifications}
To define equisingularity conditions, one can go the other way, by first considering a decomposition of a set, given by some restrictive geometric conditions on the pieces of this decomposition, and then trying to see what kind of equisingular properties have been obtained along the pieces.  
This is how the theory of stratification works. 
As a matter of fact, the theory of stratification was initiated by H. Whitney in 1957 in \cite{Whi57}, a few years before the theory of equisingularity.  
We refer to \cite{Tro20} for a complete overview on the subject, nevertheless for the reader convenience, we recall now the notions of stratification and of Whitney stratification. 

%%%%%%%%%%%ù
\begin{definition}
Let $X$ be a tame set (in $\R^n$ or $\CC^n$). A \emph{stratification} of $X$ is a decomposition of $X$ (one may also take a decomposition of the ambient space, compatible with $X$, that is, $X$ is union of some strata) into a finite number of smooth manifolds $(X_i)_{i\in \{0, \ldots, p\}}$ in the same tame class as $X$, called the \emph{strata} of the stratification, such that for any $i,j\in \{0, \ldots, p\}$, $i\not=j$, if $X_i\cap \overline{X_j}\not=\emptyset$, then $X_i\subset \overline{X_j}\smallsetminus X_j$. One then says that $(X_i,X_j)$ is a pair of \emph{adjacent strata}. 
\end{definition}

%%%%%%%%%%  
\begin{definition}\label{def. Whitney stratification}
A stratification is a \emph{Whitney stratification}, or is \emph{Whitney regular}, if any pair $(X_i,X_j)$ of adjacent strata satisfies \emph{Whitney’s condition $(b)$} at any point 
$y\in X_i$ : for all sequences $(x_i)_{i\in \N} \in X_j$ and $(y_i)_{i\in \N} \in X_i$, both having limit $y$, such that the sequence of tangent spaces $(T_{x_i}X_j)_{i\in \N}$
 tends to $T$ and the sequences of lines $(x_iy_i)_{i\in \N}$ tend to $L$ (in the correspondig Grassmannian manifold), one
has $L\subset T$.
\end{definition}

The definition of Whitney stratifications is based on a differential condition relating to the pairs of adjacent strata, a condition that can be checked by elementary calculus, at least in simple situations. It has been proved by J. Mather in \cite{Mat12} and R. Thom in \cite{Tho69} that along a stratum of one of its Whitney stratification  a set is topologically equisingular. This result, called the Thom-Mather theorem, is much more accurate, but we need not to go into details here. 
Since one knows that every real tame set admits a Whitney stratification by \cite{vdDMil,Loi98}, the Thom-Mather theorem provides another proof than the one based on the cell decomposition theorem, of the existence of topologically equisingular decompositions of real tame sets. In the same way, in the complex case, the above mentioned canonical construction of B. Teissier in \cite{Tei82} of an equisingular decomposition is the construction of a canonical minimal Whitney stratification. 
In the Section \ref{subsection. Complex numerical invariants} we will explain how B. Teissier obtained canonicity for Whitney stratifications, using algebraic characterization of Whitney regularity by numerical invariants. 

Let us notice that L\^ e D. T. and B. Teissier in \cite{LeTei83} gave a converse statement to the Thom-Mather theorem: a certain topological triviality (that we don't describe here) along the strata of a stratification implies the Whitney regularity of this stratification.

%%%%%%%%%%%%%%%%%%%%%%%%%%%%%%%%
\subsection{t-stratifications}
%%%%%%%%%%%%%%%%%%%%%%%%%%%%%%%%

The existence of Whitney stratifications for tame sets has been obtained in the $p$-adic context in \cite{CluComLoe12}, and in \cite{For17} for fields of Laurent series $k((t))$, where $k$ is an arbitrary field of characteristic $0$, with a proof similar to the classical proof for real tame sets, which consists in proving, starting from a given stratification, that the set of points in a stratum $X_i$ at which the pair of adjacent strata $(X_i,X_j)$ does not satisfy Whitney regularity, is a tame set of positive codimension in $X_i$
(see \cite{Loi98,vdDMil}). 

On the other way, non-Archimedean methods has been used to prove the existence of Whitney stratifications of real tame sets, in \cite{Hal14,BraHal22}, and even of \emph{Lipschitz stratifications} in \cite{HalYin18} (see \cite{Mos85,Par94} for the definition and the existence of Lipschitz stratifications, and \cite{Par93} or \cite[Section 4.5.3]{Par2021} for a survey on Lipschitz equisingularity). Let's now give a hint of the method, due to I. Halupczok. 

In the real or in the complex case, a Whitney stratification is first given by a differential incidence property for pairs of adjacent strata, as in Definition 
\ref{def. Whitney stratification}, and then a strong topological triviality is deduced, by the Thom-Mather theorem. The idea of I. Halupczok in \cite{Hal14,BraHal22} was to proceed in the opposite direction, by first considering  a stratification for tame non-Archimedean sets having a strong metric triviality property along its strata, and then proving that such a stratification yields a Whitney stratification when passing to the real or the complex case (see also \cite{Hal14b,Hal23} for enlightening surveys). Such a stratification is called a \emph{t-stratification} in \cite{Hal14}, and a \emph{riso-stratification} in the canonical improved version \cite{BraHal22} (in the canonical case an additional property has to be satisfied by the valued field $K$, the spherically completeness). Here is the definition.

Let $K$ be a Henselian valued field of equi-characteristic $0$, and denote $\mathrm{val}_K$ its (additive) valuation. 
For a tuple $a = (a_, \ldots, a_n)\in K^n$, we set $\mathrm{v} (a) := \min_i (\mathrm{val}_K(a_i))$. An (open) ball in $K^n$ is a set of the form
$$
B(a,\lambda) := \{ x\in K^n,  \mathrm{v}(x- a) > \lambda \},
$$
for some $a \in  K^n$ and $\lambda \in \Gamma_K$. Let us denote 
$$
\mathrm{res}: \mathcal{O}_K\to k
$$
the quotient residue map. 

%%%%%%%%%%%%%%%
\begin{definition}\label{def. t-stratification}
 With the notation above, Let $X$ be a tame set in $K^n$.
\begin{itemize}
\item[-]
A map $\Phi\colon X \to Y$ between sets $X, Y \subseteq K^n$ is called a \emph{risometry} if for every $x, x' \in X$ with $x \ne x'$, we have
$$\mathrm{v} ((\Phi(x')-\Phi(x))-(x'-x)) >\mathrm{v}(x'-x).$$

\item[-]
Let $B \subseteq K^n$ be a ball.
 Given a sub-vector space $V$ of $k^n$, we
say that $X$ is  \emph{$V$-trivial on $B$} ($\dim(V)$-translatable in \cite{Hal14} and $V$-riso trivial in \cite{BraHal22}) if there exists a lift $\widetilde{V} \subseteq K^n$ of $V$ (i.e. satisfying $\mathrm{res}(\widetilde{V}) = V$) and a risometry $\Phi: B \to B$ such that $\Phi^{-1}(X \cap B)$ is $\widetilde{V}$-translation invariant, i.e.
for every pair of points $x,y\in B$ satisfying $y-x\in \widetilde{V}$, we have 
$$ x \in \Phi^{-1}(X\cap B)\Longleftrightarrow y\in \Phi^{-1}(X\cap B).$$ 
%
%\begin{figure}[h]
%    \centering
%  \hskip1cm  
%  \includegraphics[scale=0.55]{fig1.pdf}
%  \vskip0cm
%   \caption{$1$-riso-triviality on a ball}
%    \label{f.risotriviality}
%\end{figure}
\item[-]
The \emph{triviality space of $X$ on $B$} is the maximal (with respect to inclusion) subspace $V \subseteq k^n$ such that $X$ is $V$-trivial on $B$; we denote it by $V_B(X)$. Such a maximal subspace exists by \cite[Proposition~2.3.12]{BraHal22}.
We say that $X$ is \emph{$d$-trivial} on $B$ if $V_B(X) \ge d$.

\item[-] A t-\emph{stratification of $K^n$ compatible with} $X$, is a partition $(S_i)_{i\in \{0, \ldots, n\}}$ of $K^n$ with tame $S_i$ (the strata of the stratification), such that for each $d\le n$
\begin{itemize}
\item[(1)] $\dim(S_d)=d$ or $S_d=\emptyset$,
\item[(2)] for any ball $B \subset S_d \cup \ldots \cup S_n$, the family 
$(S_d, \ldots , S_n,X)$ is $d$-trivial on $B$.
%(for the same $V$ and the same risometry).
\end{itemize}

\end{itemize}
\end{definition}

Then by \cite[Theorem 1.1]{Hal14}, and \cite[Section 4.4]{BraHal22}, for the canonical case, any tame set of $K^n$ has a t-stratification. 
Moreover by \cite[Theorem 7.11]{Hal14} and \cite[Theorem 4.6.4]{BraHal22}, starting from a tame set $X$ of $k=\R$ or $\CC$, one can associate to $X$ a non-Archimedean tame set $X^*\subset K^n$ for a valued field $K$ with residual field $k$, providing a t-stratification $(S^*_0,\ldots, S^*_n)$ of $K^n$ compatible with $X^*$. Then going back to $X$ by showing that the formula defining $(S^*_0,\ldots, S^*_n)$ can be realized in $k$, finally gives a Whitney stratification $(S_0,\ldots, S_n)$ of $k^n$ compatible with $X$.

%%%%%%%%%%%%%%%%%%%%%%%%%%%%%%%%%%%%
\section{Numerical invariants}\label{section. Numerical invariants}
%%%%%%%%%%%%%%%%%%%%%%%%%%%%%%%%%%%%
Apart from stratification theory, another possible direction to define a notion of  similarity for the germs of a tame set $X$ along one of its (tame and smooth) subset $Y$ consists in imposing the constancy along $Y$ of certain numerical invariants associated to the family of germs $(X,y)_{y\in Y}$. This have to be done with the idea that this constancy has to be deeply related to the regularity of a stratification having $X$ and $Y$ as adjacent strata.

%%%%%%%%%%%%%%%%%%%%%%
\subsection{Complex numerical invariants}\label{subsection. Complex numerical invariants}
%%%%%%%%%%%%%%%%%%%%%%
In the complex setting the most elementary invariant is the \emph{local multiplicity} $e(X,y)$ of the germ $(X,y)$. The local  multiplicity may be defined in this context by purely algebraic means, but in view to generalization to the real or the non-Archimedean cases, we rather give its geometrical definition: assume $X\subset \CC^n$, and $\dim_\CC(X)=d$,  then $e(X,y)$ is the degree of the branched covering $\pi_P: X\cap B(y,\varepsilon)\to P$, where $P$ is a generic affine space of dimension $d$ containing $y$, and $\pi_P$ the orthogonal projection on $P$, for $\varepsilon>0$ as small as desired (indeed this degree does not depend on $\varepsilon$ small enough). In other words:
\begin{definition}\label{def. multiplicity}
With the above notation, the \emph{local multiplicity} $e(X,y)$ of the complex tame germ $(X,y)$ is the number of intersection points $X\cap B(y,\varepsilon)\cap \pi_P^{-1}(z)$, for $P$ generic, $z$ generic in $P\cap B(y,\varepsilon)$, and for any $\varepsilon>0$ small enough.
\end{definition}
In \cite{Zar71}, O. Zariski noticed that ``{\sl Any definition of equisingularity should imply equimultiplicity, at the very least}", that is the constancy of $Y\ni y \mapsto e(X,y)$. He then asked a natural but ambitious question (Question A of \cite{Zar71}): 

{\sl For $(X,y)$ and $(X',y')$ two germs of complex hypersurfaces of $\CC^n$, is it true that the existence of a homeomorphism of germs 
$(\CC^n,X,y)\to (\CC^n,X',y')$ implies $e(X,y)=e(X',y')$?} 

Notice that having the same embedded topological type is necessary so that the question makes sense, since for instance $X=\{(x,y)\in \CC^2, y^2-x^3=0\}$ is homeomorphic to $X'=\CC^2$, whereas $e(X,0)=2$ and $e(X',0)=1$.

In its full generality, this question is still unsolved (for surveys on this question, one can refer to \cite{Eyr07, Eyr16, Bob2022}). Some progress has been made essentially by reinforcing the topological assumption, replacing the topological type by the Lipschitz type in Zariski question A, that is assuming that there exists a bilipschitz homeomorphism from $(\CC^n,X,y)$ to $(\CC^n,X',y')$ (Question AL$(n,d)$ of \cite{FerSam23}), or even from $(X,y)$ to $(X',y')$ (Question AL$(d)$ of \cite{FerSam23}). It has been proved in \cite{BobFerSam18} that the answer to AL$(2)$ is yes, but the answer to A$(d)$ is negative in general for $d\ge 3$ by \cite{BirFerSamVer20}.  However, by \cite{Com98}, it is known that, for arbitrary analytic sets of $\CC^n$ (and not only for hypersurface germs), if the multiplicities of $X$ and $X'$ are sufficiently closed, then A$(d)$ has a positive answer. In particular, for a Mostowski's Lipschitz stratification of a complex analytic set $X$, since along one of its stratum $Y$ one has a Lipschitz deformation of $X$ by isotopy, $Y$ is equimultiple along $X$. 
This result was in fact known since H. Hironaka has proven in \cite{Hir69} that along one of its Whitney stratification a complex analytic set is equimultiple, and since Lipschitz stratifications are also Whitney regular. 

One can refer to the survey \cite{FerSam23} for the bilipschitz version of Zariski's question on multiplicity. 

Another way to weaken Zariski's question A on multiplicity is to consider the situation of a family of hypersurface germs along a strata $Y$ of a regular stratification of $X$. More accurately, we consider an analytic mapping $f:\CC^n\times \CC\to \CC$ and set $X_y=\{z\in \CC^n, \ f(z,y)=0\}$, $Y=\{0\}\times \CC$, assuming that $f(0,y)=0$ and that $0$ is an isolated singularity of $X_y$ (at least locally at $y=0$ in $Y$). 
 Here regular means any regularity condition implying topological equisingularity of the family. One such regularity condition is that the Milnor number $\mu(X_y)$ of  $X_y$ at $0$ is constant along $Y$.
 The $\mu$-constant condition is located as follows in the hierarchy of regular conditions for hypersurfaces: 
 
 \smallskip
{\sl  Zariski's generic equisingularity $\Longrightarrow$ Whitney regularity $\Longleftrightarrow$ the constancy of Teissier's sequence $\mu^*(X_y)$ $\Longrightarrow$ the constancy of $\mu(X_y)$ $\Longrightarrow$ (embedded) topological equisingularity.  }
 \smallskip
 
The first $\Longrightarrow$ is due to \cite{Spe75}, the $\Longleftrightarrow$ is due to \cite{BriSpe76} and \cite{Tei71}, 
the second $\Longrightarrow$ comes from the fact that $\mu(X_y)$ is a term of the sequence $\mu^*(X_y)$, 
and the last  $\Longrightarrow$ is due to \cite{LeRam76} (for $n\not=3$). 
Moreover, all those implications are strict by \cite{BriSpe75}.

A remarkable positive answer to Zariski's multiplicity question has been given in the family case, by J. F. de Bobadilla: 

%%%%%%%%%%%%%
\begin{theorem}[\cite{Bob24}]
An analytic  $\mu$-constant family of hypersurfaces with isolated singularities is equimultiple. 
\end{theorem}

In general, for a stratification of a complex tame set of any codimension in $\CC^n$, a complete algebraic characterization of Whitney regularity has been given in \cite{Tei82} by B. Teissier: the equimultiplicity along the strata of all polar varieties of the set, defined as follows. 
Let $(X_i)_{i\in \{0,\ldots, p\}}$ be a stratification of a tame set $X\subset \CC^n$, and $(Y,X)$ a pair of adjacent strata of this stratification. 
The (germ of the) \emph{polar variety} $\mathcal{P}_P(X,y)$ of $X$ at $y\in Y$, associated to a (generic) projection $\pi_P:X\to P$, on the
vector spaces $P$ containing $Y$, with  $\dim(P)\in [\dim(Y)+1,n] $, is
$$
\mathcal{P}_P(X,y) := \mathrm{adh}\{x\in X, \ \dim(T_xX \cap P^\perp) \ge \dim(X)-\dim(P)+1  \},
$$   
and for $\dim(P)>\dim(X)$, $\mathcal{P}_P(X,y) :=\mathrm{adh}(X,y)$.
Let us denote by $e^*(X,y)$ the sequence of  multiplicities of $\mathcal{P}_P(X,y)$, for $\dim(P)=\dim(Y)+1, \ldots, n$. 

%%%%%%%%%%%%%
\begin{theorem}[\cite{Tei82, HenMer83, HenMerSab84}]\label{thm. Polar varieties}
With the notation above, the stratification $(X_i)_{i\in \{0,\ldots, p\}}$ is whitney regular if and only if for any adjacent pair $(Y,X)$ of this stratification, $Y\ni y\mapsto e^*(X,y) )$ is constant. 
\end{theorem}

This algebraic characterization of Whitney regularity provides the canonical minimal Whitney stratification mentioned in Section \ref{subsection. Whitney stratifications}.

%%%%%%%%%%%%%%%%%%%%%%
\subsection{Real and non-Archimedean numerical invariants}
%%%%%%%%%%%%%%%%%%%%%%
In the real case the interplay between equisingularity conditions and numerical invariants loses its simplicity (for a survey on this question one may refer to \cite{Com15}). The reason for this is twofold; on the one hand, in the real context we no longer have natural rigid invariants available, and on the other hand, even when we found convenient real substitutes to complex invariants, their regularity along the strata of a stratification does not imply in general that the stratification itself is regular.

 As an instance of this double principle, consider the question of what should be a real substitute of the local multiplicity for real a tame germ $(X,y)$. Going back to the geometric Definition \ref{def. multiplicity} of the complex local multiplicity leads to the introduction of real tame constructible functions. Given a real tame set $X\subset \R^n$, a \emph{tame constructible function} on $X$ is a function that can be written as $\displaystyle \sum_{j=1}^m \alpha_j \mathbf{1}_{Z_j}$, with $\alpha_j\in \R$, $(Z_j)_{j\in \{1, \ldots, m\}}$ a decomposition of $X$ into tame sets, and $\mathbf{1}_{Z_j}$ the characteristic function of $Z_j$. Denoting by $\mathcal{C}(X)$ the additive group of constructible functions on $X$, and by $\chi$ the Euler-Poincaré characteristic defined in any tame structure (see \cite[Chapter 4, Section 2]{vdD}), one defines a functor from the category of tame objects with tame morphisms to the category of abelian groups 
\[
\begin{tikzcd} %[column sep=huge,row sep=huge]
X \arrow[r] \arrow[d,swap,"f"] &
   \mathcal{C}(X)  \arrow[d,"f_*"] \\
W  \arrow[r] &  \mathcal{C}(W)
\end{tikzcd}
\]
where $f_*\left(\displaystyle \sum_{j=1}^m \alpha_j \mathbf{1}_{Z_j}\right)(w):=
\displaystyle \sum_{j=1}^m \alpha_j \chi (f^{-1}(w)\cap Z_j)$. 
Now, for $i\in \{0, \ldots, n\}$, for generic affine spaces $P$ of dimension $i$ containing $y$, the above diagram can be localized (see \cite[Definition 2.7]{ComMer08} or \cite[Definition 2.2.1]{Com15}) to the following diagram corresponding to tame germs 
\[
\begin{tikzcd} %[column sep=huge,row sep=huge]
(X,y) \arrow[r] \arrow[d,swap,"\pi_P"] &
   \mathcal{C}(X,y)  \arrow[d,"\pi_{P*}"] \\
(P,y)  \arrow[r] &  \mathcal{C}(P,y)
\end{tikzcd}
\]
with $\pi_{P*}(\displaystyle \sum_{j=1}^m \alpha_j \mathbf{1}_{(Z_j,y)})(w):=\displaystyle \sum_{j=1}^m \alpha_j \chi (\pi_P^{-1}(w)\cap Z_j\cap B(y,r))$, $r$ small enough, and $\vert  w  \vert \ll~r$ (so that this definition does not depend on the choices of $r$ and $w$, for $P$ generic). 
 
Then, for a given $i$-dimensional tame germ $(Z,y)$, let us denote $\Theta_i(Z,y)$ its \emph{local density}, that is the limit 
 $\Theta_i(Z,y):=\displaystyle \lim_{\varepsilon \to 0}\mathrm{Vol}_i(Z\cap B(y,\varepsilon))/\mathrm{Vol}_i(B^i(y,\varepsilon))$, and let us denote by $\theta_i$ the integral of constructible functions on $(P,y)$, $\dim(P)=i$, with respect to the local density $\Theta_i$, that is $\displaystyle  \theta_i(\displaystyle \sum_{j=1}^m \alpha_j \mathbf{1}_{Z_j})=\displaystyle \sum_{j=1}^m \alpha_j \Theta(Z_j,y)$. The existence of all those data is deeply a consequence of the tameness of the objects involved.

Finally, one can define the \emph{polar invariants $\sigma_i$} of the germ $(X,y)$, for $i=0, \ldots, n$, as 
 $$
 \sigma_i(X,y):=
 \int_{P\in G(i,n)} \theta_i(\pi_{P^*}(\mathbf{1}_{(X,y)})) \ \mathrm{d}P,
 $$ 
where $\mathrm{d}P$ is the invariant measure, under the orthogonal group, of the Grassmann manifold $G(i,n)$ of $i$-dimensional planes (containing $y$). 
\begin{remark}
 In the tame complex case, since for generic $P$ in the corresponding complex Grassmann manifold $G(d,n)$, where $d=\dim_\CC(X,y)$, 
$\pi_{P^*}(\mathbf{1}_{(X,y)})$ is generically equal to $e(X,y)$, we get 
$\sigma_d(X,y)=e(X,y)$, showing that $\sigma_d$ is a good real substitute of the local multiplicity.   
\end{remark}

On the other hand, when $(X,y)$ is a complex tame germ of $\CC^n$ of complex dimension $d$, one knows by \cite{Dra69} that $\Theta_{2d}(X,y)=e(X,y)$, and therefore $\Theta_{2d}(X,y)$ is another good real substitute of the local multiplicity.  
It turns out that even in the real context $\sigma_d=\Theta_d$, a formula called the \emph{local Cauchy-Crofton formula} (see \cite[Theorems 1.10 and 1.16]{Com00}). In the non Archimedean context, using $p$-adic or motivic integration (see \cite{CluLoe2015, Hal2018}), one can define the polar invariants $\sigma_d$ and the density $\Theta_d$ in a similar way, and they still are equal (see \cite[Theorem 6.2.1]{CluComLoe12} for the $p$-adic case, and \cite{For17} for the general non-Archimedean case). 

The interest of the sequence of polar invariants is that their variation can be studied through the geometry of the discriminants of generic projections, in the spirit of Zariski's original definition of equisingularity, and in the spirit of Hironaka's paper \cite{Hir69}, where equimultiplicity is obtained along Whitney strata through the normal pseudo-flatness of the discriminant. In the real context this study has been done, leading to the real version of the equimultiplicity statement of  \cite{Hir69}, where constancy has to be replaced by continuity, and $e(X,y)$ replaced by $\Theta_d(X,y)$.

\begin{theorem}[\cite{Com00}, real equimultiplicity]\label{thm. real equimultiplicity}
Along the strata of a Verdier stratification of a closed tame set $X$, $y\mapsto \Theta_d(X,y)$ is continuous.  
\end{theorem}

Theorem \ref{thm. real equimultiplicity} has been improved by G. Valette in \cite{Val08}, by another method based on the study of controlled vector fields along Whitney strata. It is indeed shown in \cite{Val08} that Verdier regularity may be replaced in Theorem \ref{thm. real equimultiplicity} by Whitney regularity for sets definable in polynomially bounded o-minimal structures, and that along the strata of Verdier regular stratifications, $y\mapsto \Theta_d(X,y)$ is a Lipschitz map.  
Notice that in the complex context, as shown by B. Teissier and J.-P. Henry-M. Merle in \cite{HenMer83,Tei82}, those two regularity conditions coincide, although Verdier regularity is strictly stronger than Whitney regularity in real tame geometry. Notice also that the continuity of the density function is a rather sensitive property in real geometry, since counter examples exist for Whitney stratifications in non polynomially bounded structures (see \cite{TroVal17}).  And on the other hand, the continuity of the density function along the strata of a stratification does not in general imply that this stratification is even $(a)$-regular (a regularity condition much less stronger than Whitney regularity) (see the example in the introduction of \cite{ComMer08}).  

We have noticed above that the invariant $\sigma_d$ may be defined in the complex case using complex projections, and is equal to the local multiplicity. We can define, in the same way, using complex projections, the complex version $\widetilde{\sigma}^*$ of the sequence $\sigma^*:=(\sigma_0, \ldots, \sigma_n)$. If $(X,y)$ is a germ of complex hypersurface of $\CC^n$ with isolated singularity, then the elements of  $\widetilde{\sigma}^*(X,y)$ are respectively equal, up to a constant and sign, to the elements of the sequence $\mu^*(X,y)$ of generalized Milnor numbers introduced by B. Teissier in \cite{Tei73} (see \cite[page 26]{Com15}, \cite[Section 2.2]{ComMer08}). Moreover B. Teissier has shown in \cite{Tei73} that the constancy of the sequence $\mu^*$ along the strata of a stratification of a tame complex set implies that the stratification is Whitney regular, and J. Brian\c con and J.-P. Speder has proven the converse in \cite{BriSpe76}.  

More generally, for tame set of any dimension in $\CC^n$, the sequence  $\widetilde{\sigma}^*$ of complex polar invariants have been studied first by M. Kashiwara in \cite{Kas73}, and enter in the definition of the sequence $e^*$ of the multiplicities of the polar varieties (see \cite[Theorem 6.1.9]{LeTei81},  \cite[Theorem 4.1.1]{LeTei83}, \cite[page 27]{Com15}, \cite[Section 2.2]{ComMer08} for more details). Since, as stated in Theorem \ref{thm. Polar varieties}, the constancy of the sequence $e^*$ along the strata of a stratification is equivalent to the Whitney regularity of this stratification,
the sequence $\widetilde{\sigma}^*$ is also constant along Whitney strata.

Again, we have a real version of the constancy of the sequences $e^*$ and $\sigma^*$ along Whitney strata, where constancy has to be replaced by continuity (see Theorems 4.9 and 4.10 of \cite{ComMer08}):

\begin{theorem}[real version of \cite{Tei82}]\label{thm. real polar varieties}
Along the strata of a Verdier stratification of a closed tame set $X$, $y\mapsto \sigma^*(X,y)$ is continuous.  
\end{theorem}

By \cite{NhaVal18} Verdier regularity can be relaxed to Whitney regularity in polynomially bounded o-minimal structure.  Notice also that the sequence $\sigma^*$ may be continuous along the strata of a stratification, but this stratification may not even be $(a)$-regular (see again the example in the introduction of \cite{ComMer08}). 

In \cite{Dut19}, N. Dutertre defined (see \cite[Definition 4.5]{Dut19}) a sequence $L^*$ of real invariants obtained as mean values (with respect to integration over the Grassmann manifold of vector planes of $\R^n$) of weighted densities of images of polar varieties, as real substitutes of the sequence $e^*$ of multiplicity of the complex polar varieties. Then N. Dutertre proved that the elements of $L^*$ are linear combinations of the elements of $\sigma^*$, and in view of Theorem \ref{thm. real polar varieties}, this shows that the sequence $L^*$ is continuous along Verdier strata (and in view of  \cite{NhaVal18}, continuous along Whitney strata in polynomially o-minimal structures).  

Still in the real framework, several invariants of singularities have been recently introduced, based on localization of certain curvatures, such as Lipschitz-Killing curvatures, in particular in \cite{ComMer08, Dut08, Dut15, Dut19, Dut20, Dut22}. In general such local invariants can be expressed using the invariants of the sequence $\sigma^*$. In \cite{KulLer21, BurKulLer22}
some interesting perspectives have been opened in the direction of integral geometry in the non-Archimedean framework, leaving the possibility, after the local density $\Theta_d$, of translating into this context the real invariants of singularities cited above. 
 
As already mentioned, in the complex framework,  Theorem \ref{thm. Polar varieties} provides a canonical  minimal Whitney stratification throughout the constancy of numerical invariants. A possible non-Archimedean counterpart of this point of view is achieved in \cite{BraHal22} by the construction of the Riso Tree of a tame non-Archimedean subset $X$ of $K^n$. Here is the construction of the Riso-Tree of $X$ (see \cite[Definition 3.1.5]{BraHal22} or \cite[Definition 7.1.2]{Hal23}). First of all, notice that the valuative balls $B(a,\lambda)\subset K^n$ forms a
tree $\mathcal{B}$, where each node is such a ball and the hierarchy represented by the branches connecting two nodes is the inclusion of balls. The leaves of this tree are therefore the points of $K^n$.  
%
%%%%%%%%%%%%%%%%
\begin{definition}
With the notation of Definition \ref{def. t-stratification}, the \emph{Riso-Tree} $\mathrm{Tr}(X)$ of $X$ is the partition  $\mathrm{Tr}_0,\ldots ,\mathrm{Tr}_n$ of $\mathcal{B}$ defined by
$\mathrm{Tr}_d := \{B \in \mathcal{B}; \dim(V_B(X))=d \}$.
\end{definition}
The crucial property of  $\mathrm{Tr}(X)$ is that for every $d$, there exists a $d$-dimensional tame set $Y_d \subset K^n$ such that every
ball $B \in  \mathrm{Tr}_d$ contains at least one point of $Y_d$ (see \cite[Corollary 7.5.4]{Hal23}). Then the strata of a t-stratification as in Definition \ref{def. t-stratification} can be build from the sets $Y_0, \ldots, Y_n$, in a canonical way. 

To conclude our overview, let us notice that the $0$-component $\mathrm{Tr}_0$ of $\mathrm{Tr}$ has been used in \cite{ComHal22} to define a notion of motivic number of connected components for tame non-Archimedean sets, and then to give a global motivic non-Archimedean version of the local Cauchy-Crofton formula $\sigma_d=\Theta_d$, as well as a non-Archimedean version of some real invariants called the \emph{Vitushkin variations} (see \cite[Chapter 3]{YomCom04} for the definition of the real Vitushkin invariants and \cite[Definition 5.7]{ComHal22} for the motivic case), both definitions being based on motivic integration. 

%------
% Insert acknowledgments and information
% regarding funding at the end of the last
% section, i.e., right before the bibliography.
%------

%\begin{ack}
%We thank X.
%\end{ack}

%\begin{funding}
%This work was partially supported by~\ldots
%\end{funding}

%------
% Insert the bibliography.
%------

\bibliographystyle{ems}
\bibliography{Teissier}

%\begin{thebibliography}{99}

%------ Example for a paper in journal:
% \bibitem{article1}
% A.~Petrunin, Parallel transportation for Alexandrov space with curvature bounded below.
% \emph{Geom. Funct. Anal.} \textbf{8} (1998), no.~1, 123--148.
% \Zbl{0903.53045} \MR{1601854}

%------ Example for a book:
% \bibitem{book1}
% W.~P. Ziemer, \emph{Weakly differentiable functions}.
% Grad. Texts in Math. 120,  Springer, New York, 1989.
%\Zbl{0692.46022} \MR{1014685}

%------ Example for a paper in a book:
% \bibitem{incollection1}
% J.~S. Milne, Introduction to Shimura varieties.
% In \emph{Harmonic analysis, the trace formula, and Shimura varieties},
% edited by M.~W. Marcellin and E.~Giorgi, pp. 265--378,
% Clay Math. Proc. 4, Amer. Math. Soc., Providence, RI, 2005.
% \Zbl{1148.14011} \MR{2192012}

%------ Example for a preprint on arXiv:
% \bibitem{preprint1}
% D.~V. Nguyen, S.~K. Chilappagari, M.~W. Marcellin, and B.~Vasic,
% LDPC codes from latin squares free of small trapping sets,
% 2010, \href{http://arxiv.org/abs/1008.4177}{arXiv:1008.4177}.

%------ Example for a report:
% \bibitem{report1}
% J.~Schöberl, Commuting quasi-interpolation operators.
% Technical report isc-01-10-math, Texas A\&M University, 2001,
% \url{www.isc.tamu.edu/publications-reports/tr/0110.pdf}.

%------ Example for a thesis:
% \bibitem{thesis1}
% E.~Giorgi, \emph{The geometric universe}.
% Ph.D. thesis, University of Maryland, College Park, 2002.

%\end{thebibliography}

\end{document}